\documentclass[12pt]{article}
\usepackage{latexsym}
\usepackage{amsfonts}
\usepackage{amsmath}
\usepackage{amssymb}
\usepackage{graphicx}
\usepackage{latexsym}
\usepackage{amsfonts}
\usepackage{graphicx}
\usepackage{psfrag}

\input{tcilatex}
\begin{document}

\qquad 

\thispagestyle{empty}

\begin{center}
{\Large \textbf{\ Positive Linear Maps and Perturbation bounds of Matrices }}

\vskip0.2inR. Sharma, R. Kumari\\[0pt]
Department of Mathematics \& Statistics\\[0pt]
Himachal Pradesh University\\[0pt]
Shimla -5,\\[0pt]
India - 171005\\[0pt]
email: rajesh.sharma.hpn@nic.in
\end{center}

\vskip1.5in \noindent \textbf{Abstract. \ \ \ }We show how positive unital
linear maps can be used to obtain lower bounds for the maximum distance
between the eigenvalues of two normal matrices. Some related bounds for the
spread and condition number of Hermitian matrices are also discussed here.



\bigskip

\bigskip

\bigskip

\bigskip

\bigskip

\bigskip

\bigskip

\bigskip

\bigskip

\section{\ Introduction}

\setcounter{equation}{0} \ \ \ \ Recently, Bhatia and Sharma \cite{1,4,5}
have shown that how positive unital linear maps can be used to obtain matrix
inequalities. In particular, they have obtained some old and new lower
bounds for the spread of a matrix. In this paper we show that their
technique can be extended and positive unital linear maps can also be used
to study the spectral variations of Hermitian and normal matrices.\vskip%
0.0in\noindent Let $\mathbb{M}(n)$ be the algebra of all $n\times n$ complex
matrices. Let $\left\langle x,y\right\rangle $ be the standard inner product
on $C^{n}$ defined as $\left\langle x,y\right\rangle =\sum\limits_{i=1}^{n}%
\overline{x_{i}}y_{i},$ and let $\left\Vert x\right\Vert =\left\langle
x,x\right\rangle ^{\frac{1}{2}}.$ The numerical range of an element $A\in 
\mathbb{M}(n)$ is the set 
\begin{equation*}
W\left( A\right) =\left\{ \left\langle x,Ax\right\rangle :\left\Vert
x\right\Vert =1\right\} .
\end{equation*}%
The Toeplitz-Hausdorff Theorem \cite{7,12} says that $W\left( A\right) $ is
a convex subset of the complex plane for all $A\in \mathbb{M}(n).$ For a
normal matrix $A,$%
\begin{equation*}
W\left( A\right) =Co\left( \sigma \left( A\right) \right)
\end{equation*}%
where $Co\left( \sigma \left( A\right) \right) $ denotes the convex hull of
the spectrum $\sigma \left( A\right) $ of $A$. For non-normal matrices, $%
W\left( A\right) $ may be bigger than $Co\left( \sigma \left( A\right)
\right) .$ The diameter of $W\left( A\right) $ is defined as 
\begin{equation*}
\text{diam }W\left( A\right) =\underset{i,j}{\max }\left\{ \left\vert
z_{i}-z_{j}\right\vert :z_{i},z_{j}\in W\left( A\right) \right\} .
\end{equation*}%
A linear map $\Phi :\mathbb{M}(n)\longrightarrow $ $\mathbb{M}(k)$ is called
positive if $\Phi \left( A\right) $ is positive semidefinite $\left( \text{%
psd}\right) $ whenever $A$ has that property, and unital if $\Phi \left(
I\right) =I.$ When $k=1$, such a map is called positive, unital, linear
functional and is denoted by the lower case letter $\varphi .$\vskip%
0.0in\noindent \noindent Bhatia and Davis \cite{3} have proved that if $\Phi 
$ is any positive unital linear map and the spectrum of any Hermitian matrix 
$A$ is contained in the interval $\left[ m,M\right] ,$ then%
\begin{equation}
\Phi \left( A^{2}\right) -\Phi \left( A\right) ^{2}\leq \frac{\left(
M-m\right) ^{2}}{4}.  \label{1.1}
\end{equation}%
Bhatia and Sharma \cite{4} have extended this for arbitrary matrices. One
more extension of \eqref{1.1} in the special case when $A$ is normal and $%
\varphi $ is linear functional is given in \cite{6}. They have augumented
this technique with another use of positive unital linear maps and showed
that if $\Phi _{1\text{ }}$and $\Phi _{2\text{ }}$ are positive unital
linear maps from $\mathbb{M}(n)$ to $\mathbb{M}(k)$, then for every
Hermitian matrix $A$ $\in $ $\mathbb{M}(n)$ we have%
\begin{equation}
\left\Vert \Phi _{1\text{ }}\left( A\right) -\Phi _{2\text{ }}\left(
A\right) \right\Vert \leq \text{diam }W\left( A\right)  \label{1.2}
\end{equation}%
where $\left\Vert \cdot \right\Vert $ denotes the spectral norm. Further, if 
$\varphi _{1\text{ }}$and $\varphi _{2\text{ }}$ are positive unital linear
functionals on $\mathbb{M}(n),$ then for every matrix $A$ in $\mathbb{M}(n)$%
\begin{equation}
\left\vert \mathbf{\ }\varphi _{1}\left( A\right) -\mathbf{\ }\varphi
_{2}\left( A\right) \right\vert \leq \text{diam }W\left( A\right) .
\label{1.3}
\end{equation}%
\ For more details, see \cite{5,6}. Using these inequalities they have
derived various old and new bounds for the spread of matrices. In a similar
spirit we discuss here perturbation bounds related to inequalities involving
positive linear maps.\vskip0.0in\noindent For an expository review of bounds
for the distance between the eigenvalues of two matrices $A$ and $B$ in
terms of expressions involving $\left\Vert A-B\right\Vert ,$ see \cite{2}.
In the present context an inequality of interest to us is due to Weyl $%
\left( 1911\right) $ which says that if $A$ and $B$ are Hermitian matrices,
then%
\begin{equation}
\left\Vert \func{Ei}g^{\downarrow }\left( A\right) -\func{Ei}g^{\downarrow
}\left( B\right) \right\Vert \leq \left\Vert A-B\right\Vert \leq \left\Vert 
\func{Ei}g^{\downarrow }\left( A\right) -\func{Ei}g^{\uparrow }\left(
B\right) \right\Vert  \label{1.4}
\end{equation}%
where $\func{Ei}g^{\downarrow }\left( A\right) \left( \func{Ei}g^{\uparrow
}\left( A\right) \right) $ denotes a diagonal matrix whose diagonal entries
are the eigenvalues of $A$ in decreasing $\left( \text{increasing}\right) $
order, see \cite{1,13}.\vskip0.0in\noindent \noindent For any two elements $%
A $ and $B$ of $\mathbb{M}(n),$ we define%
\begin{equation*}
s\left( W\left( A\right) ,W\left( B\right) \right) =\underset{i,j}{\max }%
\left\{ \left\vert w_{i}\left( A\right) -w_{j}\left( B\right) \right\vert
:w_{i}\left( A\right) \in W\left( A\right) ,w_{j}\left( B\right) \in W\left(
B\right) \right\} .
\end{equation*}%
Note that $s\left( W\left( A\right) ,W\left( B\right) \right) =$ diam $%
W\left( A\right) $ for $A=B.$\vskip0.0in\noindent We show that the
inequality \eqref{1.3} can be extended for two matrices $A$ and $B$ with
diam $W\left( A\right) $ replaced by $s\left( W\left( A\right) ,W\left(
B\right) \right) ,$ (see Theorem 2.1, below)$.$ In the special case when $A$
and $B$ are normal we get the lower bound for the maximum distance between
the eigenvalues of $A$ and $B$ (Corollary 2.1)$.$ Likewise, we obtain an
extension of \eqref{1.2} for two Hermitian matrices (Theorem 2.2).

\section{Main results}

\setcounter{equation}{0}\textbf{Theorem 2.1. }Let\textbf{\ }$\varphi _{i}:%
\mathbb{M}(n)\rightarrow 
\mathbb{C}
$ be positive unital linear functionals, $i=1,2.$ Let $A$ and $B$ be any two
elements of $\mathbb{M}\left( n\right) $. Then%
\begin{equation}
\left\vert \mathbf{\ }\varphi _{1}\left( A\right) -\mathbf{\ }\varphi
_{2}\left( B\right) \right\vert \leq s\left( W\left( A\right) ,W\left(
B\right) \right) .  \label{2.1}
\end{equation}%
\vskip0.0in\noindent \textbf{Proof. }If\textbf{\ }$A\in \mathbb{M}\left(
n\right) $ then every positive unital\textbf{\ }linear functional $\varphi
\left( A\right) $ can be expressed as the convex combination of $n$ complex
numbers, each of which is in the numerical range of $A,$ see \cite{5}.
Therefore, there exists complex numbers $z_{i}\left( A\right) $ in $W\left(
A\right) $ and $z_{j}\left( B\right) $ in $W\left( B\right) $ such that%
\begin{equation*}
\varphi _{1}\left( A\right) =\dsum\limits_{i=1}^{n}\alpha _{i}z_{i}\left(
A\right) \text{ and }\varphi _{2}\left( B\right)
=\dsum\limits_{j=1}^{n}\beta _{j}z_{j}\left( B\right) \text{ , }
\end{equation*}%
where $\alpha _{i}$ and $\beta _{j}$ are non- negative real numbers such
that $\dsum\limits_{i=1}^{n}\alpha _{i}=\dsum\limits_{j=1}^{n}\beta _{j}=1.$%
\vskip0.0in\noindent By the Toeplitz-Hausdorff Theorem, $\varphi _{1}\left(
A\right) \in W\left( A\right) $ and $\varphi _{2}\left( B\right) \in W\left(
B\right) $, and so \eqref{2.1} follows immediately. $\blacksquare $\vskip%
0.1in \noindent \textbf{Lemma 2.1. }Let $U$ and $V$ denote the convex hulls
of complex numbers $z_{i}\left( U\right) $ and $z_{j}\left( V\right) $
respectively, $i,j=1,2,...,n.$ Then, the inequality%
\begin{equation}
\left\vert \mathbf{\ }u-\mathbf{\ }v\right\vert \leq \underset{i,j}{\max }%
\left\{ \left\vert z_{i}\left( U\right) -z_{j}\left( V\right) \right\vert
\right\} \text{,}  \label{2.2}
\end{equation}%
holds true for all complex numbers $u\in U$ and $v\in V$.\vskip%
0.0in\noindent \textbf{Proof. }Since $u$ and $v$ are in the convex hulls of
complex numbers $z_{i}\left( U\right) $ and $z_{j}\left( V\right) $
respectively, we can write 
\begin{equation*}
u=\dsum\limits_{i=1}^{n}p_{i}z_{i}\left( U\right) \text{ and }%
v=\dsum\limits_{j=1}^{n}q_{j}z_{j}\left( V\right) ,
\end{equation*}%
where $p_{i}$ and $q_{j}$ are non-negative real numbers such that $%
\dsum\limits_{i=1}^{n}p_{i}=\dsum\limits_{j=1}^{n}q_{j}=1.$\vskip%
0.0in\noindent We therefore have 
\begin{eqnarray*}
\left\vert \mathbf{\ }u-\mathbf{\ }v\right\vert &=&\left\vert \mathbf{\ }%
\dsum\limits_{j=1}^{n}q_{j}\left( u-z_{j}\left( V\right) \right) \mathbf{\ }%
\right\vert \\
&\leq &\dsum\limits_{j=1}^{n}q_{j}\left\vert \mathbf{\ }u-z_{j}\left(
V\right) \mathbf{\ }\right\vert \\
&\leq &\underset{j}{\max }\left\vert u-z_{j}\left( V\right) \right\vert \\
&=&\underset{j}{\max }\left\vert \dsum\limits_{i=1}^{n}p_{i}\left(
z_{i}\left( U\right) -z_{j}\left( V\right) \right) \right\vert \\
&\leq &\underset{j}{\max }\dsum\limits_{i=1}^{n}p_{i}\left\vert \left(
z_{i}\left( U\right) -z_{j}\left( V\right) \right) \right\vert \\
&\leq &\underset{i,j}{\max }\left\{ \left\vert z_{i}\left( U\right)
-z_{j}\left( V\right) \right\vert \right\} .
\end{eqnarray*}%
This proves the lemma. $\blacksquare $\vskip0.1in \noindent \textbf{%
Corollary 2.1. }Under the conditions of Theorem 2.1, if $A$ and $B$ are
normal matrices, then 
\begin{equation}
\left\vert \mathbf{\ }\varphi _{1}\left( A\right) -\mathbf{\ }\varphi
_{2}\left( B\right) \right\vert \leq \underset{i,j}{\max }\left\vert \lambda
_{i}\left( A\right) -\lambda _{j}\left( B\right) \right\vert ,  \label{2.3}
\end{equation}%
where $\lambda _{i}\left( A\right) $ and $\lambda _{j}\left( B\right) $ are
the eigenvalues of $A$ and $B$, respectively.\vskip0.0in\noindent \textbf{%
Proof}. If $A$ is normal, then numerical range of $A$ is the convex polygon
spanned by the eigenvalues of $A.$ So, $W\left( A\right) $ and $W\left(
B\right) $ are the convex hulls of the eigenvalues $\lambda _{i}\left(
A\right) $ and $\lambda _{j}\left( B\right) ,$ respectively. It follows from
above lemma that 
\begin{equation*}
s\left( W\left( A\right) ,W\left( B\right) \right) =\underset{i,j}{\max }%
\left\vert \lambda _{i}\left( A\right) -\lambda _{j}\left( B\right)
\right\vert .
\end{equation*}%
The assertions of the corollary now follows from Theorem 2.1. $\blacksquare $%
\vskip0.1in \noindent \textbf{Theorem 2.2. }Let $\Phi _{1\text{ }}$and $\Phi
_{2\text{ }}$ be any two positive unital linear maps from $\mathbb{M}(n)$
into $\mathbb{M}(k)$. Let $A$ and $B$ be any two Hermitian elements of $%
\mathbb{M}(n).$ Then%
\begin{equation}
\left\Vert \Phi _{1\text{ }}\left( A\right) -\Phi _{2\text{ }}\left(
B\right) \right\Vert \leq \left\Vert \func{Ei}g^{\downarrow }\left( A\right)
-\func{Ei}g^{\uparrow }\left( B\right) \right\Vert .  \label{2.4}
\end{equation}%
\vskip0.0in\noindent \textbf{Proof}. If $A$ is an $n\times n$ Hermitian
matrix then $\lambda _{n}^{\downarrow }\left( A\right) I\leq A\leq \lambda
_{1}^{\downarrow }\left( A\right) I.$ The linear map $\Phi _{1\text{ }}$
preserves order and take the identity $I$ in $\mathbb{M}(n)$ to $I$ in $%
\mathbb{M}(k).$ So we have $\lambda _{n}^{\downarrow }\left( A\right) I\leq
\Phi _{1\text{ }}\left( A\right) \leq \lambda _{1}^{\downarrow }\left(
A\right) I.$ Likewise, we have $\lambda _{1}^{\uparrow }\left( B\right)
I\leq \Phi _{2\text{ }}\left( B\right) \leq \lambda _{n}^{\uparrow }\left(
B\right) I.$ It then follows that%
\begin{equation*}
\lambda _{n}^{\downarrow }\left( A\right) -\lambda _{n}^{\uparrow }\left(
B\right) \leq \Phi _{1\text{ }}\left( A\right) -\Phi _{2\text{ }}\left(
B\right) \leq \lambda _{1}^{\downarrow }\left( A\right) -\lambda
_{1}^{\uparrow }\left( B\right) .
\end{equation*}%
Therefore,%
\begin{equation*}
-k\leq \Phi _{1\text{ }}\left( A\right) -\Phi _{2\text{ }}\left( B\right)
\leq k
\end{equation*}%
where 
\begin{equation*}
k=\max \left\{ \left\vert \lambda _{n}^{\downarrow }\left( A\right) -\lambda
_{n}^{\uparrow }\left( B\right) \right\vert ,\left\vert \lambda
_{1}^{\downarrow }\left( A\right) -\lambda _{1}^{\uparrow }\left( B\right)
\right\vert \right\}
\end{equation*}%
Further, if $X$ is Hermitian and $\pm X\leq kI$ then $\left\Vert
X\right\Vert \leq k$, and therefore%
\begin{equation*}
\left\Vert \Phi _{1\text{ }}\left( A\right) -\Phi _{2\text{ }}\left(
B\right) \right\Vert \leq k.
\end{equation*}%
The assertions of the theorem now follow from the fact that%
\begin{equation*}
\left\Vert \func{Ei}g^{\downarrow }\left( A\right) -\func{Ei}g^{\uparrow
}\left( B\right) \right\Vert =\underset{j}{\max }\left\vert \lambda
_{j}^{\downarrow }\left( A\right) -\lambda _{j}^{\uparrow }\left( B\right)
\right\vert =k.\text{ }\blacksquare
\end{equation*}%
\vskip0.0in\noindent We note that the inequality \eqref{2.4} and the second
inequality \eqref{1.4} are independent. The maps 
\begin{equation*}
\Phi _{1\text{ }}\left( A\right) =\frac{1}{n-1}\left( trA-A\right) \text{
and }\Phi _{2\text{ }}\left( B\right) =B
\end{equation*}%
are positive unital linear maps. For these maps, the inequality \eqref{2.4}
becomes 
\begin{equation}
\left\Vert \func{Ei}g^{\downarrow }\left( A\right) -\func{Ei}g^{\uparrow
}\left( B\right) \right\Vert \geq \frac{1}{n-1}\left\Vert A-B+n\left( B-%
\frac{trA}{n}\right) \right\Vert .  \label{2.5}
\end{equation}%
For $A=B$, the inequality \eqref{2.5} gives 
\begin{equation}
\left\Vert \func{Ei}g^{\downarrow }\left( A\right) -\func{Ei}g^{\uparrow
}\left( B\right) \right\Vert \geq \frac{n}{n-1}\left\Vert A-\frac{trA}{n}%
\right\Vert ,  \label{2.6}
\end{equation}%
while Weyl's inequality \eqref{1.4} gives $\left\Vert \func{Ei}g^{\downarrow
}\left( A\right) -\func{Ei}g^{\uparrow }\left( B\right) \right\Vert \geq 0.$%
\vskip0.0in\noindent But for $B=\frac{trA}{n}I,$ we respectively have from %
\eqref{2.5} and \eqref{1.4}, 
\begin{equation*}
\left\Vert \func{Ei}g^{\downarrow }\left( A\right) -\func{Ei}g^{\uparrow
}\left( B\right) \right\Vert \geq \frac{1}{n-1}\left\Vert A-B\right\Vert
\end{equation*}%
and 
\begin{equation*}
\left\Vert \func{Ei}g^{\downarrow }\left( A\right) -\func{Ei}g^{\uparrow
}\left( B\right) \right\Vert \geq \left\Vert A-B\right\Vert .
\end{equation*}%
\vskip0.1in \noindent \noindent \noindent Choosing different linear maps in
Theorem 2.2 and Corollary 2.1, we can obtain various interesting
inequalities which provide lower bounds for the maximum distance between
eigenvalues of two normal matrices. We demonstrate some special cases here.%
\vskip0.2in\noindent \noindent Choose $\varphi _{1}\left( A\right) =a_{ii}$
and $\varphi _{2}\left( B\right) =b_{jj}$ in \eqref{2.3}$,$ we have%
\begin{equation}
\underset{i,j}{\max }\left\vert \lambda _{i}\left( A\right) -\lambda
_{j}\left( B\right) \right\vert \geq \max_{i,j}\left\vert
a_{ii}-b_{jj}\right\vert .  \label{2.7}
\end{equation}%
Let $D$ be the diagonal part of $A.$ From \eqref{2.7} we have%
\begin{equation}
\max_{i,j}\left\vert \lambda _{i}\left( A\right) -\lambda _{j}\left(
D\right) \right\vert \geq \max_{i,j}\left\vert a_{ii}-a_{jj}\right\vert .
\label{2.8}
\end{equation}%
Note that \eqref{2.8} provides a refinement of the inequality spd$\left(
A\right) \geq \underset{i,j}{\max }\left\vert a_{ii}-a_{jj}\right\vert .$ By
using spectral theorem, $a_{ii}=\dsum\limits_{i=1}^{n}\lambda _{i}\left(
A\right) p_{i}$ where $p_{i}$ are non-negative real numbers such that $%
\dsum\limits_{i=1}^{n}p_{i}=1,$ therefore%
\begin{eqnarray*}
\left\vert \lambda _{i}\left( A\right) -a_{jj}\right\vert &=&\left\vert
\dsum\limits_{j=1}^{n}p_{j}\left( \lambda _{i}\left( A\right) -\lambda
_{j}\left( A\right) \right) \right\vert \leq
\dsum\limits_{j=1}^{n}p_{j}\left\vert \left( \lambda _{i}\left( A\right)
-\lambda _{j}\left( A\right) \right) \right\vert \\
&\leq &\underset{i,j}{\max }\left\vert \lambda _{i}\left( A\right) -\lambda
_{j}\left( A\right) \right\vert =\text{spd}\left( A\right) .
\end{eqnarray*}%
So spd$\left( A\right) \geq \underset{i,j}{\max }\left\vert \lambda
_{i}\left( A\right) -\lambda _{j}\left( D\right) \right\vert \geq \underset{%
i,j}{\max }\left\vert a_{ii}-a_{jj}\right\vert .$\vskip0.0in\noindent Let $%
A=D+N$ . If $N$ is also a normal matrix, as in case of circulant and
Hermitian matrices, then%
\begin{equation*}
\underset{i,j}{\max }\left\vert \lambda _{i}\left( A\right) -\lambda
_{j}\left( N\right) \right\vert \geq \underset{i}{\max }\left\vert
a_{ii}\right\vert .
\end{equation*}%
Let $B=\frac{A+A^{\ast }}{2}$ and $C=\frac{A-A^{\ast }}{2i}.$ If $A$ is
normal, $\func{Re}\left( \lambda _{i}\left( A\right) \right) =\lambda
_{i}\left( B\right) $ and $\func{Im}\left( \lambda _{i}\left( A\right)
\right) =\lambda _{i}\left( C\right) .$We therefore have%
\begin{equation*}
\underset{i,j}{\max }\left\vert \lambda _{i}\left( A\right) -\func{Re}\left(
\lambda _{j}\left( A\right) \right) \right\vert \geq \underset{i,j}{\max }%
\left\vert a_{ii}-\func{Re}a_{jj}\right\vert ,
\end{equation*}%
\begin{equation*}
\underset{i,j}{\max }\left\vert \lambda _{i}\left( A\right) -\func{Im}\left(
\lambda _{j}\left( A\right) \right) \right\vert \geq \underset{i,j}{\max }%
\left\vert a_{ii}-\func{Im}a_{jj}\right\vert .
\end{equation*}%
For arbitrary matrices, we have%
\begin{equation*}
\underset{i,j}{\max }\left\vert \lambda _{i}\left( \frac{A+A^{\ast }}{2}%
\right) -\lambda _{j}\left( \frac{A-A^{\ast }}{2i}\right) \right\vert \geq 
\underset{i,j}{\max }\left\vert \func{Re}a_{ii}-\func{Im}a_{jj}\right\vert
\end{equation*}%
Note that if $A$ is normal, $\lambda _{i}\left( \frac{A+A^{\ast }}{2}\right)
-\lambda _{j}\left( \frac{A-A^{\ast }}{2i}\right) =\func{Re}\left( \lambda
_{i}\left( A\right) \right) -\func{Im}\left( \lambda _{j}\left( A\right)
\right) .$\noindent \vskip0.1in \noindent We now obtain some more
inequalities in the following corollaries.\noindent \vskip0.1in \noindent 
\textbf{Corollary 2.2.} Let $A=\left[ a_{ij}\right] $ and $B=\left[ b_{ij}%
\right] $ be Hermitian matrices. Then%
\begin{equation}
\left\Vert \func{Ei}g^{\downarrow }\left( A\right) -\func{Ei}g^{\uparrow
}\left( B\right) \right\Vert \geq \frac{1}{2}\left\vert \alpha +\beta \pm 
\sqrt{\left( \alpha -\beta \right) ^{2}+4\left\vert a_{ij}+b_{ij}\right\vert
^{2}}\right\vert  \label{2.9}
\end{equation}%
where $\alpha =a_{ii}-b_{jj}$\textbf{, }$\beta =a_{jj}-b_{ii}$ and $i\neq j.$%
\vskip0.0in\noindent \textbf{Proof}. The maps 
\begin{equation*}
\Phi _{1\text{ }}\left( A\right) =\left[ 
\begin{array}{cc}
a_{ii} & a_{ij} \\ 
a_{ji} & a_{jj}%
\end{array}%
\right] \text{ \ \ \ and \ \ \ }\Phi _{2\text{ }}\left( B\right) =\left[ 
\begin{array}{cc}
b_{jj} & -b_{ij} \\ 
-b_{ji} & b_{ii}%
\end{array}%
\right] .
\end{equation*}%
are positive unital linear maps, and 
\begin{equation*}
\Phi _{1\text{ }}\left( A\right) -\Phi _{2\text{ }}\left( B\right) =\left[ 
\begin{array}{cc}
a_{ii}-b_{jj} & a_{ij}+b_{ij} \\ 
a_{ji}+b_{ji} & a_{jj}-b_{ii}%
\end{array}%
\right] .
\end{equation*}%
is a Hermitian matrix with eigenvalues 
\begin{equation*}
\frac{1}{2}\left( \alpha +\beta \pm \sqrt{\left( \alpha -\beta \right)
^{2}+4\left\vert a_{ij}+b_{ij}\right\vert ^{2}}\right) \text{ \ }.
\end{equation*}%
So $\left\Vert \Phi _{1\text{ }}\left( A\right) -\Phi _{2\text{ }}\left(
B\right) \right\Vert $\bigskip $=\frac{1}{2}\left\vert \alpha +\beta \pm 
\sqrt{\left( \alpha -\beta \right) ^{2}+4\left\vert a_{ij}+b_{ij}\right\vert
^{2}}\right\vert .$ The inequality \eqref{2.9} now follows from Theorem 2.2. 
$\blacksquare $\vskip0.0in\noindent In the special case when $A=B$,
inequality \eqref{2.9} gives Mirsky bound \cite{11} for the spread of $A$.
It also follows from \eqref{2.9} that 
\begin{equation*}
\left\Vert \func{Ei}g^{\downarrow }\left( A\right) -\func{Ei}g^{\uparrow
}\left( D\right) \right\Vert \geq \sqrt{\left( a_{ii}-a_{jj}\right)
^{2}+\left\vert a_{ij}\right\vert ^{2}}
\end{equation*}%
and%
\begin{equation*}
\left\Vert \func{Ei}g^{\downarrow }\left( A\right) -\func{Ei}g^{\uparrow
}\left( N\right) \right\Vert \geq \frac{1}{2}\left\vert a_{ii}+a_{jj}\pm 
\sqrt{\left( a_{ii}-a_{jj}\right) ^{2}+16\left\vert a_{ij}\right\vert ^{2}}%
\right\vert .
\end{equation*}%
\vskip0.1in \noindent \textbf{Corollary 2.3. }Let $A$ and $B$ be normal
matrices. Let $I$ and $J$ be any two subsets of $\left\{ 1,2,...,n\right\} $
and let $\left\vert I\right\vert $ and $\left\vert J\right\vert $ denote the
cardinality of $I$ and $J$ . Then%
\begin{equation}
\underset{i,j}{\max }\left\vert \lambda _{i}\left( A\right) -\lambda
_{j}\left( B\right) \right\vert \geq \max_{i,j}\left\vert \frac{1}{%
\left\vert I\right\vert }\dsum\limits_{i,j\in I}a_{ij}-\frac{1}{\left\vert
J\right\vert }\dsum\limits_{i,j\in J}b_{ij}\right\vert .  \label{2.10}
\end{equation}%
\vskip0.0in\noindent \noindent \textbf{Proof}. Choose 
\begin{equation*}
\varphi _{1}\left( A\right) =\frac{1}{\left\vert I\right\vert }%
\dsum\limits_{i,j\in I}a_{ij}\text{ and }\varphi _{2}\left( B\right) =\frac{1%
}{\left\vert J\right\vert }\dsum\limits_{i,j\in J}b_{ij}
\end{equation*}%
and use \eqref{2.3}, we immediately get \eqref{2.10}. $\blacksquare $\vskip%
0.0in\noindent A special case of Corollary 2.3 when $A$ $=$ $B$ is Theorem
2.2 of Johnson et al \cite{8}$.$ The inequality \eqref{2.10} for $B=D$ gives%
\begin{equation*}
\underset{i,j}{\max }\left\vert \lambda _{i}\left( A\right) -\lambda
_{j}\left( D\right) \right\vert \geq \max_{i,j}\left\vert \frac{1}{%
\left\vert I\right\vert }\dsum\limits_{i,j\in I}a_{ij}-\frac{1}{\left\vert
J\right\vert }\dsum\limits_{i\in J}a_{ii}\right\vert .
\end{equation*}%
It also follows from \eqref{2.10} that%
\begin{equation*}
\underset{i,j}{\max }\left\vert \lambda _{i}\left( A\right) -\func{Re}\left(
\lambda _{j}\left( A\right) \right) \right\vert \geq \max_{i,j}\left\vert 
\frac{1}{\left\vert I\right\vert }\dsum\limits_{i,j\in I}a_{ij}-\frac{1}{%
\left\vert J\right\vert }\dsum\limits_{i,j\in J}\func{Re}a_{ij}\right\vert ,
\end{equation*}%
\begin{equation*}
\underset{i,j}{\max }\left\vert \lambda _{i}\left( A\right) -\func{Im}\left(
\lambda _{j}\left( A\right) \right) \right\vert \geq \max_{i,j}\left\vert 
\frac{1}{\left\vert I\right\vert }\dsum\limits_{i,j\in I}a_{ij}-\frac{1}{%
\left\vert J\right\vert }\dsum\limits_{i,j\in J}\func{Im}a_{ij}\right\vert
\end{equation*}%
and%
\begin{equation*}
\underset{i,j}{\max }\left\vert \func{Re}\left( \lambda _{i}\left( A\right)
\right) -\func{Im}\left( \lambda _{j}\left( A\right) \right) \right\vert
\geq \max_{i,j}\left\vert \frac{1}{\left\vert I\right\vert }%
\dsum\limits_{i,j\in I}\func{Re}a_{ij}-\frac{1}{\left\vert J\right\vert }%
\dsum\limits_{i,j\in J}\func{Im}a_{ij}\right\vert .
\end{equation*}%
\vskip0.1in \noindent \textbf{Corollary 2.4. }Let $A$ and $B$ be normal
matrices. Then%
\begin{equation}
\underset{i,j}{\max }\left\vert \lambda _{i}\left( A\right) -\lambda
_{j}\left( B\right) \right\vert \geq \left\vert \frac{1}{n-1}\sum_{i\neq
j}a_{ij}+\frac{1}{n}\sum_{i,j}\left( b_{ij}-a_{ij}\right) \right\vert .
\label{2.11}
\end{equation}%
\noindent \vskip0.0in\noindent \noindent \textbf{Proof}. For the positive
unital linear functionals 
\begin{equation*}
\varphi _{1\text{ }}\left( A\right) =\frac{1}{n}\left( trA-\frac{1}{n-1}%
\sum_{i\neq j}a_{ij}\right) ,\text{ \ }\varphi _{2}\left( B\right) =\frac{1}{%
n}\sum_{i,j}b_{ij}\text{ \ \ and }\varphi _{3}\left( A\right) \text{ }=\frac{%
1}{n}\sum_{i,j}a_{ij}\text{,}
\end{equation*}%
we have%
\begin{eqnarray*}
\left\vert \varphi _{1\text{ }}\left( A\right) -\varphi _{2}\left( B\right)
\right\vert &=&\left\vert \varphi _{1\text{ }}\left( A\right) -\varphi _{3%
\text{ }}\left( A\right) +\varphi _{3\text{ }}\left( A\right) -\varphi
_{2}\left( B\right) \right\vert \\
&=&\left\vert \frac{1}{n-1}\sum_{i\neq j}a_{ij}+\frac{1}{n}\sum_{i,j}\left(
b_{ij}-a_{ij}\right) \right\vert .
\end{eqnarray*}%
The assertions of the corollary now follows from the inequality \eqref{2.3}. 
$\blacksquare $\vskip0.0in\noindent \noindent Theorem 5 of Merikoski and
Kumar \cite{9} is a special case of our Corollary 2.4, $A=B$. For $A=D$ and $%
B=A,$%
\begin{equation*}
\max_{i,j}\left\vert \lambda _{i}\left( A\right) -\lambda _{j}\left(
D\right) \right\vert \geq \frac{1}{n}\left\vert \sum_{i\neq
j}a_{ij}\right\vert .
\end{equation*}%
and%
\begin{equation*}
\underset{i,j}{\max }\left\vert \func{Re}\left( \lambda _{i}\left( A\right)
\right) -\func{Im}\left( \lambda _{j}\left( A\right) \right) \right\vert
\geq \left\vert \frac{1}{n-1}\sum_{i\neq j}\func{Re}a_{ij}+\frac{1}{n}%
\sum_{i,j}\left( \func{Im}a_{ij}-\func{Re}a_{ij}\right) \right\vert .
\end{equation*}%
\vskip0.1in \noindent \textbf{Corollary 2.5.} Let $A$ and $B$ be normal
matrices. Then%
\begin{equation}
\underset{i,j}{\max }\left\vert \lambda _{i}\left( A\right) -\lambda
_{j}\left( B\right) \right\vert \geq \frac{1}{2}\left\vert \left( \alpha
+\beta \right) +\left( a_{ij}e^{i\theta }+a_{ji}e^{-i\theta }\right)
\right\vert .  \label{2.12}
\end{equation}%
\noindent \vskip0.0in\noindent \noindent \textbf{Proof}. Let 
\begin{equation*}
\varphi _{1\text{ }}\left( A\right) =\frac{1}{2}\left(
a_{ii}+a_{jj}+a_{ij}e^{i\theta }+a_{ji}e^{-i\theta }\right) \text{ \ \ and \
\ }\varphi _{2}\left( B\right) =\frac{1}{2}\left( b_{ii}+b_{jj}\right) .
\end{equation*}%
$\varphi _{1\text{ }}\left( A\right) $ and $\varphi _{2\text{ }}\left(
B\right) $ are positive unital linear functionals, see \cite{5}$.$ The
inequality \eqref{2.12} follows from \eqref{2.3}. $\ \blacksquare $\vskip%
0.0in\noindent Let $B=D.$ Then, \eqref{2.12} gives%
\begin{equation*}
\max_{i,j}\left\vert \lambda _{i}\left( A\right) -\lambda _{j}\left(
D\right) \right\vert \geq \frac{1}{2}\underset{i\neq j}{\max }\left\vert
a_{ij}e^{i\theta }+a_{ji}e^{-i\theta }\right\vert .
\end{equation*}%
So,%
\begin{equation*}
\max_{i,j}\left\vert \lambda _{i}\left( A\right) -\lambda _{j}\left(
D\right) \right\vert \geq \frac{1}{2}\underset{i\neq j}{\max }\left(
\left\vert a_{ij}\right\vert +\left\vert a_{ji}\right\vert \right) .
\end{equation*}%
If $A$ is Hermitian, we have%
\begin{equation*}
\max_{i,j}\left\vert \lambda _{i}\left( A\right) -\lambda _{j}\left(
D\right) \right\vert \geq \underset{i\neq j}{\max }\left\vert
a_{ij}\right\vert .
\end{equation*}%
Likewise, we can see that

\begin{equation*}
\max_{i,j}\left\vert \lambda _{i}\left( A\right) -\lambda _{j}\left(
D\right) \right\vert \geq \underset{p\neq q}{\max }\left\vert \frac{1}{n}%
\dsum\limits_{i,j}a_{ij}-\frac{a_{pp}+a_{qq}}{2}\right\vert .
\end{equation*}

\section{Bounds for spread}

\setcounter{equation}{0}The spread of $A$, denoted by spd$\left( A\right) ,$
is defined as 
\begin{equation*}
\text{spd}\left( A\right) =\underset{1\leq i,j\leq n}{\max }\left\vert
\lambda _{i}\left( A\right) -\lambda _{j}\left( A\right) \right\vert ,
\end{equation*}%
where $\lambda _{1}\left( A\right) ,...,\lambda _{n}\left( A\right) $ are
the eigenvalues of $A.$ Begining with Mirsky \cite{10} several authors have
worked on the bounds for spread of matrices, see \cite{3,4,5,6} and
references therein. We mention here some lower bounds for the spread related
to perturbation bounds. It is clear from $\left( 2.6\right) $ that for any
Hermitian element $A\in \mathbb{M}(n),$ we have%
\begin{equation*}
\text{spd}\left( A\right) \geq \frac{n}{n-1}\left\Vert A-\frac{trA}{n}%
\right\Vert 
\end{equation*}%
We prove that this inequality also holds for normal matrices.\vskip0.1in
\noindent \textbf{Theorem 3.1.} \noindent For any normal matrix $A$, we have%
\begin{equation}
\text{spd}\left( A\right) \geq \frac{n}{n-1}\left\Vert A-\frac{trA}{n}%
\right\Vert .  \label{3.1}
\end{equation}%
\vskip0.0in\noindent \textbf{Proof}. It is immediate that 
\begin{equation*}
\left\vert \lambda _{j}-\frac{trA}{n}\right\vert \leq \frac{1}{n}%
\dsum\limits_{i=1,i\neq j}^{n}\left\vert \lambda _{i}-\lambda
_{j}\right\vert \leq \frac{n-1}{n}\max_{i}\left\vert \lambda _{i}-\lambda
_{j}\right\vert ,
\end{equation*}%
for all $j=1,2...,n.$ So 
\begin{equation}
\text{spd}\left( A\right) \geq \frac{n}{n-1}\max_{j}\left\vert \lambda _{j}-%
\frac{trA}{n}\right\vert .  \label{3.2}
\end{equation}%
For a normal matrix $A$, we have%
\begin{equation}
\max_{j}\left\vert \lambda _{j}-\frac{trA}{n}\right\vert =\left\Vert A-\frac{%
trA}{n}\right\Vert .  \label{3.3}
\end{equation}%
Combining \eqref{3.2} and \eqref{3.3}$,$we immediately get \eqref{3.1}. $\
\blacksquare $\vskip0.0in\noindent It may be noted here that for a normal
matrix $A$, 
\begin{equation*}
\max_{j}\left\vert \lambda _{j}-\frac{trA}{n}\right\vert \geq \left\Vert
\varphi \left( A-\frac{trA}{n}\right) \right\Vert .
\end{equation*}%
We therefore also have 
\begin{equation}
\text{spd}\left( A\right) \geq \frac{n}{n-1}\left\Vert \varphi \left(
A\right) -\frac{trA}{n}\right\Vert .  \label{3.4}
\end{equation}%
Choose $\varphi \left( A\right) =\frac{1}{n}\dsum\limits_{i,j=1}^{n}a_{ij}$, %
\eqref{3.4} gives%
\begin{equation*}
\text{spd}\left( A\right) \geq \frac{1}{n-1}\left\vert \dsum\limits_{i\neq
j}a_{ij}\right\vert .
\end{equation*}%
This is Theorem $2.1$ of Johnson et al \cite{8} and Theorem $5$ of Merikoski
and Kumar \cite{9}$.$ We now prove a refinement of Theorem $5$ in \cite{9}$.$%
\vskip0.1in \noindent \textbf{Theorem 3.2.} Let $A$ be a normal matrix. Then%
\begin{equation}
\text{spd}\left( A\right) \geq \max_{i,j}\left\vert \lambda _{i}\left(
A\right) -\frac{1}{n-1}\dsum\limits_{k\neq j,k=1}^{n}\lambda _{k}\left(
A\right) \right\vert \geq \frac{1}{n-1}\left\vert \dsum\limits_{i\neq
j}a_{ij}\right\vert .  \label{3.5}
\end{equation}%
\vskip0.0in\noindent \noindent \textbf{Proof}. To prove first inequality %
\eqref{3.5}, note that%
\begin{eqnarray*}
\left\vert \lambda _{i}\left( A\right) -\frac{1}{n-1}\dsum\limits_{k\neq
j,k=1}^{n}\lambda _{k}\left( A\right) \right\vert  &=&\frac{1}{n-1}%
\left\vert \dsum\limits_{k\neq j,k=1}^{n}\lambda _{i}\left( A\right)
-\lambda _{k}\left( A\right) \right\vert  \\
&\leq &\frac{1}{n-1}\dsum\limits_{k\neq j,k=1}^{n}\left\vert \lambda
_{i}\left( A\right) -\lambda _{k}\left( A\right) \right\vert  \\
&\leq &\max_{j}\left\vert \lambda _{i}\left( A\right) -\lambda _{j}\left(
A\right) \right\vert .
\end{eqnarray*}%
For $B=\frac{1}{n-1}\left( trA-A\right) $ and $\varphi _{1}=\mathbf{\ }%
\varphi _{2},$ the inequality \eqref{2.3} gives%
\begin{equation}
\underset{i,j}{\max }\left\vert \lambda _{i}\left( A\right) -\lambda
_{j}\left( B\right) \right\vert \geq \frac{n}{n-1}\left\vert \mathbf{\ }%
\varphi _{1}\left( A\right) -\mathbf{\ }\frac{trA}{n}\right\vert .
\label{3.6}
\end{equation}%
The eigenvalues of $B$ are $\frac{1}{n-1}\dsum\limits_{k\neq
j,k=1}^{n}\lambda _{k}\left( A\right) ,$ $j=1,2,...,n.$ Choose $\varphi
_{1}\left( A\right) =\frac{1}{n}\dsum\limits_{i.j}a_{ij}$ in \eqref{3.6}, we
immediately get the second inequality in \eqref{3.5}. $\ \blacksquare $\vskip%
0.1in \noindent Analogous bound for the ratio spread $\left( \frac{\lambda
_{\max }}{\lambda _{\min }}\right) $of a positive definite matrix is proved
in the following theorem.\vskip0.1in \noindent \textbf{Theorem 3.3. }Let $%
\Phi :\mathbb{M}(n)\rightarrow \mathbb{M}(k)$ be a positive unital linear
map. Let $A\in \mathbb{M}(n)$ be a positive definite matrix. For $m\leq
A\leq M$, we have%
\begin{equation}
\left( \frac{m}{M}\right) ^{\frac{n-1}{n}}\leq \left( \det A\right) ^{\frac{%
-1}{n}}\Phi \left( A\right) \leq \left( \frac{M}{m}\right) ^{\frac{n-1}{n}}.
\label{3.7}
\end{equation}%
\vskip0.0in\noindent \textbf{Proof}. Let $\lambda _{i}$ be the eigenvalues
of $A.$ It is clear that $m\leq \lambda _{i}\leq M,$ $i=1,2,...,n.$
Therefore,%
\begin{equation*}
m^{n-1}M\leq \lambda _{1}\lambda _{2}...\lambda _{n}\leq M^{n-1}m\text{ .}
\end{equation*}%
So,%
\begin{equation*}
m^{\frac{n-1}{n}}M^{\frac{1}{n}}\leq \left( \det A\right) ^{\frac{1}{n}}\leq
M^{\frac{n-1}{n}}m^{\frac{1}{n}}
\end{equation*}%
and therefore%
\begin{equation}
\frac{1}{M^{\frac{n-1}{n}}m^{\frac{1}{n}}}\leq \frac{1}{\left( \det A\right)
^{\frac{1}{n}}}\leq \frac{1}{m^{\frac{n-1}{n}}M^{\frac{1}{n}}}.  \label{3.8}
\end{equation}%
Also, $m\leq A\leq M$ therefore%
\begin{equation}
m\leq \Phi \left( A\right) \leq M.  \label{3.9}
\end{equation}%
The inequality \eqref{3.7} now follows from \eqref{3.8} and \eqref{3.9}. $\
\blacksquare $\vskip0.1in \noindent \noindent Let $\Phi :\mathbb{M}%
(n)\longrightarrow $ $\mathbb{M}(k)$ be a positive unital linear map. Let $A$
be a Hermitian element of $\mathbb{M}(n)$ such that $mI\leq A\leq MI$.
Bhatia and Davis \cite{3} have proved that 
\begin{equation}
\Phi \left( A^{2}\right) -\Phi \left( A\right) ^{2}\leq \left( M-\Phi \left(
A\right) \right) \left( \Phi \left( A\right) -m\right) .  \label{3.10}
\end{equation}%
We use \eqref{3.10}, and obtain a refinement of \eqref{1.1} for positive
definite matrices under a certain condition.\vskip0.1in \noindent \textbf{%
Theorem 3.4. }Let $\Phi :\mathbb{M}(n)\longrightarrow $ $\mathbb{M}(k)$ be a
positive unital linear map and $A$ be any positive semidefinite matrix in $%
\mathbb{M}(n)$ such that $mI\leq A\leq MI$. If $\Phi \left( A^{2}\right)
\geq 2\Phi \left( A\right) ^{2}$ and $\Phi \left( A\right) >0,$ then%
\begin{equation}
\Phi \left( A^{2}\right) \leq \left( M-m\right) \Phi \left( A\right) .
\label{3.11}
\end{equation}%
\vskip0.0in\noindent \noindent \textbf{Proof}. From \eqref{3.10}, we have%
\begin{equation}
\Phi \left( A^{2}\right) \leq \left( m+M-mM\Phi \left( A\right) ^{-1}\right)
\Phi \left( A\right) .  \label{3.12}
\end{equation}%
If $\Phi \left( A^{2}\right) \geq 2\Phi \left( A\right) ^{2},$ then%
\begin{equation}
\Phi \left( A\right) \leq \sqrt{\Phi \left( A^{2}\right) -\Phi \left(
A\right) ^{2}}.  \label{3.13}
\end{equation}%
From \eqref{1.1} and \eqref{3.13}, we get that 
\begin{equation*}
\Phi \left( A\right) \leq \frac{M-m}{2}.
\end{equation*}%
Therefore%
\begin{equation*}
m+M-mM\Phi \left( A\right) ^{-1}\leq m+M-\frac{2mM}{M-m}\leq M-m.
\end{equation*}%
Since $\Phi \left( A\right) $ and $\Phi \left( A\right) ^{-1}$ commute,
therefore%
\begin{equation}
\left( m+M-mM\Phi \left( A\right) ^{-1}\right) \Phi \left( A\right) \leq
\left( M-m\right) \Phi \left( A\right) .  \label{3.14}
\end{equation}%
From \eqref{3.12} and \eqref{3.14}, we get \eqref{3.11}. \ $\blacksquare $%
\vskip0.1in We now show that above theorem provides a refinement of the
inequality \eqref{1.1} for linear functionals.\vskip0.1in \noindent \textbf{%
Corollary 3.1.} Let $\varphi :\mathbb{M}(n)\longrightarrow $ $\mathbb{C}$ be
a positive unital linear functional and $A$ be any positive semidefinite
matrix in $\mathbb{M}(n)$ such that $mI\leq A\leq MI$. If $\varphi \left(
A^{2}\right) \geq 2\varphi \left( A\right) ^{2}$ and $\varphi \left(
A\right) >0,$ then%
\begin{equation}
\varphi \left( A^{2}\right) -\varphi \left( A\right) ^{2}\leq \left( \frac{%
\varphi \left( A^{2}\right) }{2\varphi \left( A\right) }\right) ^{2}\leq 
\frac{\left( M-m\right) ^{2}}{4}.  \label{3.15}
\end{equation}%
\vskip0.0in\noindent \textbf{Proof}. The second inequality in \eqref{3.15}
follows from inequality \eqref{3.11}. The first inequality holds if and only
if%
\begin{equation*}
4\varphi \left( A^{2}\right) \varphi \left( A\right) ^{2}-4\varphi \left(
A\right) ^{4}\leq \varphi \left( A^{2}\right) ^{2}
\end{equation*}%
if and only if%
\begin{equation*}
\left( \varphi \left( A^{2}\right) -2\varphi \left( A\right) ^{2}\right)
^{2}\geq 0
\end{equation*}%
This is true. \ $\blacksquare $\vskip0.0in\textbf{Example}. Let%
\begin{equation*}
A=\left[ 
\begin{array}{ccc}
2 & 2 & 1 \\ 
2 & 2 & 1 \\ 
1 & 1 & 1%
\end{array}%
\right] ,\text{ \ \ }A^{2}=\left[ 
\begin{array}{ccc}
9 & 9 & 5 \\ 
9 & 9 & 5 \\ 
5 & 5 & 3%
\end{array}%
\right] 
\end{equation*}%
Let $\varphi \left( A\right) =a_{ii}$, then $\varphi \left( A^{2}\right)
\geq 2\varphi \left( A\right) ^{2}$. From \eqref{1.1} we have $M-m\geq 4.4721
$. From\eqref{3.15}, $M-m\geq 4.5$. Note that $M-m\geq 4.5616.$ \vskip%
0.2in\noindent \textbf{Acknowledgements}. The authors are grateful to Prof.
Rajendra Bhatia for the useful discussions and suggestions, and I.S.I. Delhi
for a visit in January 2015 when this work had begun. The support of the
UGC-SAP is also acknowledged.

\end{document}